\documentclass[12pt]{article}
\usepackage{amsthm,amsmath,amssymb}
\newtheorem{theorem}{Theorem}[section]

\theoremstyle{definition}

\theoremstyle{remark}

\def\R{\hbox{{\rm I}\kern-0.2em{\rm R}\kern0.2em}}
\def\D{\hbox{{\rm I}\kern-0.2em{\rm D}\kern0.2em}}

\def\be{\begin{equation}}
\def\ee{\end{equation}}

\def\({\left(}
\def\){\right)}
\def\[{\left[}
\def\]{\right]}
\def\bc{\begin{center}}
\def\ec{\end{center}}

\begin{document}

{\bf\Large \textbf{Classification of scalar third order ordinary differential equations linearizable via generalized contact transformations}}
\begin{center}

\textbf{Hina M. Dutt$^a$, Asghar Qadir$^a$} \end{center}
$^a$School of Natural Sciences\\
National University of Sciences and Technology\\
Campus H-12, 44000, Islamabad, Pakistan\\
E-mail: hinadutt@yahoo.com; aqadirmath@yahoo.com

{\bf Abstract}.

Whereas Lie had linearized scalar second order ordinary differential
equations (ODEs) by point transformations and later Chern had
extended to the third order by using contact transformation, till
recently no work had been done for higher order or systems of ODEs.
Lie had found a unique class defined by the number of infinitesimal
symmetry generators but the more general ODEs were not so
classified. Recently classifications of higher order and systems of
ODEs were provided. In this paper we relate contact
symmetries of scalar ODEs with point symmetries of reduced systems.
We define new type of transformations that build up this relation
and obtain equivalence classes of scalar third order ODEs linearizable via these
transformations. Four equivalence classes of such equations
are seen to exist.

\section{Introduction}
Symmetries of ODEs have acquired considerable interest over the last
few decades. The earlier work had dealt with second order linear
ODEs. The number of point symmetries a scalar second order ODE can
have is exactly one of $0,1,2,3,$ or $8$ \cite{sur1}. Lie
\cite{liei7,lie1,liei2,lie2} discovered that if a second order
scalar ODE possesses $8$ point symmetries, then it can be converted
to a linear equation (\emph{linearizable}) via invertible point
transformations. He proved that a scalar second order linearizable
ODE is at most cubic in the first derivative. The coefficients of
linearizable ODEs satisfy four constraints. These constraints were
reduced to two by Tress$\grave{e}$ \cite{tre1}.

Mahomed and Leach \cite{mahomedleachequiv} showed that all
linearizable scalar ODEs of order $n$ with $(n\geq3)$ can be put
into three equivalence classes with $n+1$, $n+2$ and $n+4$ Lie point
symmetries. Chern \cite{chern1,chern2} used contact transformations
to reduce scalar third order ODEs to $y'''=0$ and $y'''+y=0$. Grebot
\cite{grebot,g8} addressed the linearization problem of scalar third
order ODEs via restricted point transformations. Neut and Petitot
\cite{neut}, and independently Ibragimov and Meleshko (IM)
\cite{ibr12,ibrmel1} used the original Lie procedure of point
transformations to linearize scalar third order ODEs to the form
$y'''(x)+k_0(x)y(x)=0$. They determined the linearizability criteria
as well as the procedure for the construction of linearizing
transformation for scalar third order ODEs. IM
also used contact transformations to linearize scalar third order
ODEs. The linearization problem for scalar fourth order ODEs gets
more complicated and was dealt with by Ibragimov, Meleshko and
Suksern (IMS) \cite{ibrmelsuk,Suksern}. They obtained the linearization
criteria both by point and contact transformations. Meleshko \cite{meleshko} also provided a simple algorithm for reducing third order autonomous ODEs to Lie linearizable second order ODEs and was extended to scalar fourth order autonomous ODEs \cite{dutt}. The class of equations linearizable by Meleshko algorithm are not included in IM and IMS classes.

As such there are three equivalence classes for scalar third order
ODEs linearizable via point transformations but no classification for
these equations linearizable via contact transformations was
provided. Indeed IM \cite{ibrmel1} obtained the
linearizable form via contact transformations for a scalar third
order ODE but they neither discussed its contact symmetries nor
provided the classification for it. The subject of the present paper
is to investigate the contact symmetries of
linearizable scalar ODEs. In fact we reduce the scalar third order ODE to a system of two second order ODEs and relate contact
symmetries of the scalar ODE with point symmetries of
the reduced system. In doing so we define a new class of
transformations that are more general than contact transformations and perform group classification of scalar
third order ODEs linearizable via these transformations.

The plan of the paper is as follows. In the subsequent section we
provide some preliminaries. In the third section we define
generalized contact transformations. The canonical form of
linearizble scalar third order ODEs is obtained and group
classification is performed by using this form. The last section contains a summary discussion of
the results.

\section{Preliminaries}
A differential equation (DE) is said to be invariant
(\emph{symmetrical}) under a \emph{point transformation}
\begin{eqnarray}
t=\varphi(x,y),\quad u=\psi(x,y),\nonumber
\end{eqnarray}
if it is form invariant under that transformation.
Instead of using finite transformations we can use infinitesimal symmetry generators
\begin{eqnarray}
\textbf{X}=\xi(x,y){\partial \over \partial x}+\eta(x,y){\partial \over \partial y}.\nonumber
\end{eqnarray}
To apply the above symmetry generator on $n^{th}$ order ODEs, we
need to extend or prolong it. The prolonged generator up to order
$n$ is given by
{\small
\begin{eqnarray}
\textbf{X}^{(n)}=\xi(x,y){\partial \over \partial x}+\eta(x,y){\partial \over \partial y}+\eta^{(1)}(x,y,y'){\partial \over \partial y'}+\eta^{(2)}(x,y,y',y''){\partial \over \partial y''}\nonumber\\
+\ldots+\eta^{(n)}(x,y,y',\ldots, y^{(n)}){\partial \over \partial y^{(n)}},\nonumber
\end{eqnarray}
}
where $\eta^{(n)}$ is the extension coefficient of $\textbf{X}^{(n)}$ and is given by
\begin{eqnarray}
\eta^{(n)}={d\eta^{(n-1)}\over dx}-y^{(n)}{d\xi\over dx}.\nonumber
\end{eqnarray}
If we have a system of $m$ ODEs of order $n$ i.e., $m$ dependent variables $\textbf{y}=y_i,~(i=1,2,\ldots,m)$, then
{\small
\begin{eqnarray}
\textbf{X}^{(n)}=\xi(x,y_i){\partial \over \partial x}+\eta_i(x,y_i){\partial \over \partial y_i}+\eta^{(1)}_i(x,y_i,y'_i){\partial \over \partial y'_i}+\eta^{(2)}_i(x,y_i,y'_i,y''_i){\partial \over \partial y''_i}\nonumber\\
+\ldots+\eta_i^{(n)}(x,y_i,y'_i,\ldots, y^{(n)}_i){\partial \over \partial y^{(n)}_i}.
\end{eqnarray}
}
Now $\textbf{X}^{(n)}$ will be the symmetry generator or simply the symmetry of the system of $m$ ODEs of order $n$
\begin{eqnarray}
y_i^{(n)}=f_i(x,y_i;y_i',y_i'',\ldots,y_i^{(n-1)}),\quad (i=1,2,\ldots,m),\label{6.01}
\end{eqnarray}
if it satisfies the relations
\begin{eqnarray}
\textbf{X}^{(n)}[y_i^{(n)}-f_i(x,y_i;y_i',y_i'',\ldots,y_i^{(n-1)})]=0,\quad (i=1,2,\ldots,m),\nonumber
\end{eqnarray}
which imply that
\begin{eqnarray}
\eta_i^{(n)}=\textbf{X}^{(n-1)}f_i,\quad
(i=1,2,\ldots,m),\label{2.001}
\end{eqnarray}
which are called symmetry conditions for the system of ODEs
(\ref{6.01}).

Two equations are said to be \emph{equivalent} if there exists a
point transformation which maps one equation into the other. All
first order scalar ODEs are equivalent. However this does not hold
for second or higher order ODEs. Nevertheless all linear scalar
second order ODEs are equivalent \cite{sur} and can be transformed
to $y''=0$. For scalar $n^{th}$ order linear ODEs with $n\geq3$ we
have a theorem due to Laguerre \cite{lag1,lag2} which says that any
linear homogeneous $n^{th}$ order scalar ODE
\begin{eqnarray}
y^{(n)}+\sum_{i=0}^{n-1}k_i(x)y^{(i)}=0~,\quad n\geq3,\label{6.002}
\end{eqnarray}
can be transformed by a point transformation to an equation of the form
\begin{eqnarray}
y^{(n)}+\sum_{i=0}^{n-3}k_i(x)y^{(i)}=0.\label{6.003}
\end{eqnarray}
 Eq (\ref{6.003}) is called the \emph{Laguerre canonical form} of the linear scalar ODEs (\ref{6.002}). Hence there are three equivalence classes arising from the Laguerre canonical form for the scalar third order ODEs. These three classes are: $y'''(x)=0,$ $y'''(x)+k_0y(x)=0$ and $y'''(x)+k_0(x)y(x)=0,$ which have $7$, $5$ and $4$ point symmetries respectively.

We need to generalize contact
transformations. Hence we shall first review the definition of
these transformations. A transformation
\begin{eqnarray}
t=\varphi(x,y,z),\quad u=\psi(x,y,z),\quad v=\omega(x,y,z),\label{6.1}
\end{eqnarray}
where $z\equiv y'$, is called a \emph{contact transformation} if it
satisfies the contact condition, $v\equiv u'={du\over dt}$.

\section{Generalized contact transformations}
Consider an $n^{th}$ order scalar ODE with $(n\geq3)$
\begin{eqnarray}
y^{(n)}=f(x,y,y',y'',\ldots,y^{(n-1)}).\label{6.012}
\end{eqnarray}
We substitute $y'=z$ in the above ODE. This reduces the above scalar ODE to the following system of two ODEs of order $n-1$
\begin{eqnarray}
&&y^{(n-1)}=z^{(n-2)},\nonumber\\
&&z^{(n-1)}=f(x,y,z;z',z'',\ldots,z^{(n-2)}).\label{6.001}
\end{eqnarray}
A point transformation
\begin{eqnarray}
t=\varphi(x,y,z),\quad u=\psi(x,y,z),\quad v=\omega(x,y,z),\label{6.013}
\end{eqnarray}
 for the system (\ref{6.001}) corresponds to a \emph{generalized contact transformation} for the scalar $n^{th}$ order ODE (\ref{6.012}) with $z=y'$. The transformation (\ref{6.013}) is actually the contact transformation without the contact
 condition.

The system of equations (\ref{6.001}) is linearizable if there is a point transformation (\ref{6.013}) which reduces the system (\ref{6.001}) to a linear system of ODEs.

\subsection{Group classification}
Consider the general form of a linear scalar third order ODE
\begin{eqnarray}
y'''=\delta(x)+\sigma(x)y+\alpha(x)y'+\beta(x)y''.\label{6.3}
\end{eqnarray}
Defining $y'=z$ we reduce the order to $2$ and double the dimensions, i.e.
\begin{eqnarray}
&&y''=z',\nonumber\\
&&z''=\delta(x)+\sigma(x)y+\alpha(x)z+\beta(x)z'.\label{6.5}
\end{eqnarray}
Here the variable $z$ is actually the derivative of the variable $y$. For the purpose of the group classification we replace the variable $z$ by $y'$ and compare the corresponding system of equations (\ref{6.5}) with the second canonical form of the linear system of second order ODEs \cite{wm} given by
\begin{eqnarray}
y''=k_1(x)y'+k_2(x)z',\nonumber\\
z''=k_3(x)y'+k_4(x)z'.\label{6.10}
\end{eqnarray}
Since (\ref{6.10}) only depends explicitly on $y'$ and $z'$ so it gives $\delta(x)=0$ and $\sigma(x)=0$.
The system (\ref{6.5}) is now of the form
\begin{eqnarray}
&&y''=z',\nonumber\\
&&z''=\alpha(x)y'+\beta(x)z'.\label{2.114}
\end{eqnarray}
This is the reduced form of the systems of second order ODEs that is obtained from a scalar linear, third order ODE. We now perform the group classification on the above system. For this we suppose
{\small
\begin{eqnarray}
\textbf{X}^{(2)}=\xi(x,y,z){\partial\over \partial x}+\eta_1(x,y,z)
{\partial\over \partial y}+\eta_2(x,y,z){\partial\over \partial
z}+\eta^{(1)}_1(x,y,z){\partial\over \partial
y'}\nonumber\\
+\eta^{(1)}_2(x,y,z){\partial\over \partial
z'}+\eta^{(2)}_1(x,y,z){\partial\over \partial
y''}+\eta^{(2)}_2(x,y,z){\partial\over \partial
z''},\label{6.15}
\end{eqnarray}
}
be the symmetry generator for (\ref{2.114}), then symmetry conditions (\ref{2.001}) read as
\begin{eqnarray}
&&\eta^{(2)}_1=\textbf{X}^{(1)}(z'),\nonumber\\
&&\eta^{(2)}_2=\textbf{X}^{(1)}(\alpha(x)y+\beta(x)z'),\label{6.16}
\end{eqnarray}
where $\eta^{(2)}_1$ and $\eta^{(2)}_2$ are second
extension coefficients of the generator (\ref{6.15}) and are given by
{\small
\begin{eqnarray}
\eta^{(2)}_1=\eta_{1,xx}+y'(2\eta_{1,xy}-\xi_{,xx})+2z'\eta_{1,xz}+y''(\eta_{1,y}-2\xi_{,x}-3y'\xi_{,y}-2z'\xi_{,z})\nonumber\\
+z''(\eta_{1,z}-y'\xi_{,z})+y'^2(\eta_{1,yy}-2\xi_{,xy}-2z'\xi_{,yz})+2y'z'(\eta_{1,yz}-\xi_{,xz})\nonumber\\
+z'^2(\eta_{1,zz}-y'\xi_{,zz})-y'^3\xi_{,yy}~,\label{6.23}\\
\eta^{(2)}_2=\eta_{2,xx}+z'(2\eta_{1,xz}-\xi_{,xx})+2y'\eta_{2,xy}+z''(\eta_{2,z}-2\xi_{,x}-2y'\xi_{,y}-3z'\xi_{,z})\nonumber\\
+y''(\eta_{2,y}-z'\xi_{,y})+z'^2(\eta_{2,zz}-2\xi_{,xz}-2z'\xi_{,yz})+2y'z'(\eta_{2,yz}-\xi_{,xy})\nonumber\\
+y'^2(\eta_{2,zz}-z'\xi_{,yy})-z'^3\xi_{,zz}.\label{6.24}
\end{eqnarray}
}
The symmetry conditions (\ref{6.16}) give the following system of determining PDEs
{\footnotesize
\begin{eqnarray}
\xi_{,yy}=\xi_{,yz}=\xi_{,zz}=0,\quad \eta_{1,zz}-\xi_{,z}=0,\quad \eta_{1,xx}-\eta_{2,x}=0, \label{6.17}\\
\eta_{1,yy}-2\xi_{,xy}-\alpha\xi_{,z}=0,\quad \eta_{2,yy}-\alpha\xi_{,y}=0,\quad 2\eta_{2,zz}-2\xi_{,xz}-\xi_{,y}-2\beta\xi_{,z}=0,\label{6.18}\\
2\eta_{1,yz}-2\xi_{,xz}-2\xi_{,y}-\beta\xi_{,z}=0,\quad 2\eta_{2,yz}-2\xi_{,xy}-\beta\xi_{,y}-2\alpha\xi_{,z}=0,\label{6.19}\\
2\eta_{1,xy}-\eta_{2,y}-\xi_{,xx}+\alpha\eta_{1,z}=0,\quad 2\eta_{2,xy}+\alpha\eta_{2,z}-\beta\eta_{2,y}-\alpha\eta_{1,y}-\xi \alpha_{,x}-\alpha\xi_{,x}=0,\label{6.20}\\
2\eta_{1,xz}+\eta_{1,y}-\eta_{2,z}-\xi_{,x}+\beta\eta_{1,z}=0,\quad 2\eta_{2,xz}-\alpha\eta_{1,z}+\eta_{2,y}-\xi_{,xx}-\xi \beta_{,x}-\beta \xi_{,x}=0,\label{6.21}\\
\eta_{2,xx}-\alpha\eta_{1,x}-\beta\eta_{2,x}=0.\label{6.22}
\end{eqnarray}
}
The system of PDEs (\ref{6.17}) gives the following solution
{\small
\begin{eqnarray}
&&\xi=ya_1(x)+za_2(x)+a_3(x),\label{6.25}\\
&&\eta_1={1\over 2}z^2a_{2}+za_4(x,y)+a_5(x,y),\label{6.26}\\
&&\eta_2={1\over 2}z^2a_{2,x}+za_{4,x}+a_{5,x}+a_6(y,z),\label{6.27}
\end{eqnarray} }
where $a_i,~(i=1,2,\ldots,6)$ are arbitrary functions of their arguments.

We now assume $\beta(x)$ to be zero, nonzero constant and an arbitrary function of $x$ and consider the following cases.

\textbf{\emph{Case I}}~$\beta(x)=0$\\
The system of PDEs (\ref{6.18})$-$(\ref{6.22}) in this case takes the form
{\footnotesize
\begin{eqnarray}
\eta_{1,yy}-2\xi_{,xy}-\alpha\xi_{,z}=0,\quad \eta_{2,yy}-\alpha\xi_{,y}=0,\quad 2\eta_{2,zz}-2\xi_{,xz}-\xi_{,y}=0,\label{6.39}\\
2\eta_{1,yz}-2\xi_{,xz}-2\xi_{,y}=0,\quad 2\eta_{2,yz}-2\xi_{,xy}-2\alpha\xi_{,z}=0,\label{6.40}\\
2\eta_{1,xy}-\eta_{2,y}-\xi_{,xx}+\alpha\eta_{1,z}=0,\quad 2\eta_{2,xy}+\alpha\eta_{2,z}-\alpha\eta_{1,y}-\xi \alpha_{,x}-\alpha\xi_{,x}=0,\label{6.41}\\
2\eta_{1,xz}+\eta_{1,y}-\eta_{2,z}-\xi_{,x}=0,\quad 2\eta_{2,xz}-\alpha\eta_{1,z}+\eta_{2,y}-\xi_{,xx}=0,\label{6.42}\\
\eta_{2,xx}-\alpha\eta_{1,x}=0.\label{6.43}
\end{eqnarray}
}
The above system of PDEs is solved for different values of $\alpha(x)$.

\textbf{\emph{Case I.1}}~\emph{Both} $\alpha$ \emph{and} $\beta$ \emph{are zero}\\
In this case the system of PDEs (\ref{6.39})$-$(\ref{6.43}) reduces to
{\small
\begin{eqnarray}
\eta_{1,yy}-2\xi_{,xy}=0,\quad \eta_{2,yy}=0,\quad 2\eta_{2,zz}-2\xi_{,xz}-\xi_{,y}=0,\quad 2\eta_{1,yz}-2\xi_{,xz}-2\xi_{,y}=0,\nonumber\\ 2\eta_{2,yz}-2\xi_{,xy}=0,\quad 2\eta_{1,xy}-\eta_{2,y}-\xi_{,xx}=0,\quad 2\eta_{2,xy}=0,\nonumber\\
2\eta_{1,xz}+\eta_{1,y}-\eta_{2,z}-\xi_{,x}=0,\quad
2\eta_{2,xz}+\eta_{2,y}-\xi_{,xx}=0,\quad
\eta_{2,xx}=0.\nonumber
\end{eqnarray}
}
Solving the above system yields the following
$15$ Lie point symmetries:
{\footnotesize
\begin{eqnarray}
&&\textbf{X}_1={\partial\over \partial x},\quad
\textbf{X}_2={\partial\over \partial y},\quad \textbf{X}_3={\partial\over \partial z},\label{6.38}\\
&&\textbf{X}_{4}=x{\partial\over \partial y},\quad \textbf{X}_5=z{\partial\over \partial y},\quad \textbf{X}_6
={1\over 2}x^2{\partial\over \partial y}+x{\partial\over \partial z},\\
&&\textbf{X}_7=z{\partial\over
\partial x}+{1\over 2}z^2{\partial\over \partial y},\quad \textbf{X}_8=x{\partial\over \partial x}+{1\over 2}xz{\partial\over \partial y},\quad \textbf{X}_9= {1\over 2}xz{\partial\over \partial y}+z{\partial\over \partial z},\\
&&\textbf{X}_{10}=({1\over 2} xz+y){\partial\over \partial y},\quad \textbf{X}_{11}={1\over 2}x^2 {\partial\over \partial x}+xy {\partial\over \partial y}+y {\partial\over \partial z},\\
&&\textbf{X}_{12}=x^2{\partial\over \partial x}+{1\over 2}(x^2z+xy) {\partial\over \partial y}+xz {\partial\over \partial z},\\
&&\textbf{X}_{13}={-1\over 2}xz{\partial\over \partial x}+{1\over 4}z(xz+2y){\partial\over \partial y}+{1\over 2}z^2{\partial\over \partial z},\\
&&\textbf{X}_{14}=({1\over 2}xz+y){\partial\over \partial x}-{1\over 4}z(xz-2y){\partial\over \partial y},\\
&&\textbf{X}_{15}=({1\over 2}x^2z-xy){\partial\over \partial x}+({1\over 4}x^2z^2-y^2){\partial\over \partial y}+({1\over 2}xz^2-yz){\partial\over \partial z}.
\end{eqnarray}
}
\textbf{\emph{Case I.2}}~$\alpha=\alpha_0\neq0$, $\beta=0$\\
The system of PDEs (\ref{6.39})$-$(\ref{6.43}) in this case yields a $15-$dimensional Lie algebra. The first three operators are $\textbf{X}_1, \textbf{X}_2, \textbf{X}_3$, given by (\ref{6.38}), while the remaining $12$ operators are
{\footnotesize
\begin{eqnarray}
&&\textbf{Y}_{1}=y{\partial\over \partial y}+z{\partial\over \partial z},\label{6.43.1}\\
&&\textbf{Y}_{2}=z{\partial\over \partial y}+\alpha_0y{\partial\over \partial z},\label{6.44}\quad \textbf{Y}_{3}=e^{\sqrt{\alpha_0}x}{\partial\over \partial y}+\sqrt{\alpha_0}e^{\sqrt{\alpha_0}x}{\partial\over \partial z},\\
&&\textbf{Y}_{4}=y{\partial\over \partial x}+zy{\partial\over \partial y}+{1\over 2}(y^2\alpha_0+z^2){\partial\over \partial z},\\ &&\textbf{Y}_{5}=z{\partial\over \partial x}+{1\over 2}(y^2\alpha_0+z^2){\partial\over \partial y}+\alpha_0zy{\partial\over \partial z},\\
&&\textbf{Y}_{6}=e^{-\sqrt{\alpha_0}x}{\partial\over \partial y}-\sqrt{\alpha_0}e^{-\sqrt{\alpha_0}x}{\partial\over \partial z},\\
&&\textbf{Y}_{7}=e^{\sqrt{\alpha_0}x}{\partial\over \partial x}+\sqrt{\alpha_0}e^{\sqrt{\alpha_0}x}y{\partial\over \partial y}+\alpha_0e^{\sqrt{\alpha_0}x}y{\partial\over \partial z},\\
&&\textbf{Y}_{8}=e^{-\sqrt{\alpha_0}x}{\partial\over \partial x}-\sqrt{\alpha_0}e^{-\sqrt{\alpha_0}x}y{\partial\over \partial y}+\alpha_0e^{-\sqrt{\alpha_0}x}y{\partial\over \partial z},\\
&&\textbf{Y}_{9}={(z\sqrt{\alpha_0}-\alpha_0y)e^{\sqrt{\alpha_0}x}\over \sqrt{\alpha_0}}{\partial\over \partial y}+(z\sqrt{\alpha_0}-\alpha_0y)e^{\sqrt{\alpha_0}x}{\partial\over \partial z},\\
&&\textbf{Y}_{10}={(z\sqrt{\alpha_0}+\alpha_0y)e^{-\sqrt{\alpha_0}x}\over \sqrt{\alpha_0}}{\partial\over \partial y}-(z\sqrt{\alpha_0}-\alpha_0y)e^{-\sqrt{\alpha_0}x}{\partial\over \partial z},\\
&&\textbf{Y}_{11}={(\sqrt{\alpha_0}y-z)e^{\sqrt{\alpha_0}x}\over \sqrt{\alpha_0}}{\partial\over \partial x}+{(\alpha_0y^2-z^2)e^{\sqrt{\alpha_0}x}\over 2\sqrt{\alpha_0}}{\partial\over \partial y}+{(\alpha_0y^2-z^2)e^{\sqrt{\alpha_0}x}\over 2}{\partial\over \partial z},\\
&&\textbf{Y}_{12}={(\sqrt{\alpha_0}y+z)e^{-\sqrt{\alpha_0}x}\over \sqrt{\alpha_0}}{\partial\over \partial x}+{(\alpha_0y^2-z^2)e^{-\sqrt{\alpha_0}x}\over 2\sqrt{\alpha_0}}{\partial\over \partial y}+{(\alpha_0y^2-z^2)e^{-\sqrt{\alpha_0}x}\over 2}{\partial\over \partial z}.
\end{eqnarray}
}
\textbf{\emph{Case I.3.1}} $\alpha=(x\pm c)^m,~m\neq 2,$ or $e^x,~\beta=0$\\
This case produces a $5-$dimensional Lie algebra. The first three operators are $\textbf{X}_2, \textbf{X}_3, \textbf{Y}_1$ given by (\ref{6.38}) and (\ref{6.43.1}).

\textbf{\emph{Case I.3.2}} $\alpha=(x\pm c)^{-2},~\beta=0$\\
In this case we obtain a $6-$dimensional Lie algebra.

\textbf{\emph{Case II}} $\beta(x)\neq 0$\\
The following subcases arise:

\textbf{\emph{Case II.1}} $\alpha=0,~\beta=\beta_0\neq0$\\
The system of PDEs (\ref{6.18})$-$(\ref{6.22}) simplifies to
{\footnotesize
\begin{eqnarray}
\eta_{1,yy}-2\xi_{,xy}=0,\quad \eta_{2,yy}=0,\quad 2\eta_{2,zz}-2\xi_{,xz}-\xi_{,y}-2\beta_0\xi_{,z}=0,\nonumber\\
2\eta_{1,yz}-2\xi_{,xz}-2\xi_{,y}-\beta_0\xi_{,z}=0,\quad 2\eta_{2,yz}-2\xi_{,xy}-\beta_0\xi_{,y}=0,\nonumber\\
2\eta_{1,xy}-\eta_{2,y}-\xi_{,xx}=0,\quad 2\eta_{2,xy}-\beta_0\eta_{2,y}=0,\quad \eta_{2,xx}-\beta_0\eta_{2,x}=0,\nonumber\\
2\eta_{1,xz}+\eta_{1,y}-\eta_{2,z}-\xi_{,x}+\beta_0\eta_{1,z}=0,\quad 2\eta_{2,xz}+\eta_{2,y}-\xi_{,xx}-\beta_0 \xi_{,x}=0,\nonumber
\end{eqnarray} }
which produces a $7-$dimensional Lie algebra. The first four of these operators are $\textbf{X}_1,\textbf{X}_2,\textbf{X}_3,\textbf{Y}_1$ given by (\ref{6.38}) and (\ref{6.43.1}) while the remaining three are
{\small
\begin{eqnarray}
\textbf{Y}_{2}=x{\partial\over \partial y},\quad \textbf{Y}_{3}=(-\beta_0y+z){\partial\over \partial y},\quad \textbf{Y}_{4}={e^{\beta_0x}\over \beta_0}{\partial\over \partial y}+e^{\beta_0x}{\partial\over \partial z}.\nonumber
\end{eqnarray}
}
~\\
\textbf{\emph{Case II.2.1}}~$\alpha=\alpha_0=\beta_0\neq 0$\\
In this subcase the system of PDEs (\ref{6.18})$-$(\ref{6.22}) gives the following set of solution
{\small
\begin{eqnarray}
&&\xi=c_1,\quad \eta_1=c_2y+c_3z+c_4+c_5e^{\alpha_1x}+c_6e^{\alpha_2x},\nonumber\\
&&\eta_2=c_5\alpha_1e^{\alpha_1x}+c_6\alpha_2e^{\alpha_2x}+c_3(y+z)\alpha_0+c_2z+c_7,\nonumber
\end{eqnarray} }
where
{\small
\begin{eqnarray}
\alpha_1={1\over 2}(\alpha_0+\sqrt{\alpha_0^2+4\alpha_0}),\quad \alpha_2={1\over 2}(\alpha_0-\sqrt{\alpha_0^2+4\alpha_0})\nonumber,
\end{eqnarray}
}
and $c_i,~(i=1,2,\ldots,7)$ are arbitrary constants. This yields a $7-$dimensional Lie algebra with $\textbf{X}_1,~\textbf{X}_2,~\textbf{X}_3$ and $\textbf{Y}_1$ given by (\ref{6.38}) and (\ref{6.43.1}). The extra three operators are
{\footnotesize
\begin{eqnarray}
\textbf{Y}_2=z{\partial\over \partial y}+\alpha_0(y+z){\partial\over \partial z},\quad \textbf{Y}_3=e^{\alpha_1x}{\partial\over \partial y}+\alpha_1e^{\alpha_1x}{\partial\over \partial z},\quad \textbf{Y}_4=e^{\alpha_2x}{\partial\over \partial y}+\alpha_2e^{\alpha_2x}{\partial\over \partial z}.\nonumber
\end{eqnarray}
}
~\\
\textbf{\emph{Case II.2.2}}~$\alpha=\alpha_0\neq0$ and $\beta=\beta_0\neq0$ with $\alpha_0\neq \beta_0$\\
This case produces the following set of solution for the PDEs (\ref{6.18})$-$(\ref{6.22})
{\small
\begin{eqnarray}
&&\xi=c_1,\quad \eta_1=c_2y+c_3z+c_4+c_5e^{\beta_1x}+c_6e^{\beta_2x},\nonumber\\
&&\eta_2=c_5\beta_1e^{\beta_1x}+c_6\beta_2e^{\beta_2x}+c_3\alpha_0y+c_3\beta_0z+c_2z+c_7,\nonumber
\end{eqnarray}
}
where
{\footnotesize
\begin{eqnarray}
\beta_1={1\over 2}(\beta_0+\sqrt{\beta_0^2+4\alpha_0}),\quad \text{and}\quad \beta_2={1\over 2}(\beta_0-\sqrt{\beta^2+4\alpha_0}).\nonumber
\end{eqnarray}
}
From above we get a $7-$dimensional Lie algebras. The first four operators are given by (\ref{6.38}) and (\ref{6.43.1}).\\
~\\
\textbf{\emph{Case II.3}}~$\beta=\beta_0,~\alpha(x)=cx^m,(x+c)^m,m=1,2~$\\
This case produces a $5-$dimensional Lie algebra with first three generators $\textbf{X}_2,\textbf{X}_3,\textbf{Y}_1$.
~\\
\textbf{\emph{Case II.4.1}}~$\alpha(x)=(cx\pm d)^m,~\beta(x)=(cx\pm d)^m,~m=1,2$\\
Here we find a algebra with $6$ Lie point symmetries.\\
~\\
\textbf{\emph{Case II.4.2}}~$\alpha(x)=\beta(x)=x^{(-1)},x^{(-2)}$\\
This case lies in the above case.

The system of PDEs (\ref{6.17})$-$(\ref{6.22}) provides us four equivalence classes with $5$, $6$, $7$ and $15$ symmetries. These Lie point symmetries correspond to generalized contact symmetries for the scalar ODE (\ref{6.3}). Thus we have the following theorem.
\begin{theorem}
If a scalar third order ODE is linearizable via generalized contact transformation $(\ref{6.013})$ then it has one of $5$, $6$, $7$ and $15$ generators of generalized contact transformations.
\end{theorem}

\section{Conclusions}
Though Leach and Mahomed \cite{mahomedleachequiv} had shown that there are three equivalence classes of third order scalar ODEs linearizable via point transformations, no work on symmetry group classification of these equations linearizable via contact transformations was done. In fact IM \cite{ibrmel1} got the necessary form of a scalar third order ODE linearizable via contact transformations but there was no attempt to address the classification problem with their methods. In this paper we have found a connection between (generalized) contact transformations of system of order $n$ with point transformations of system of order $n-1$. By defining the first derivative of dependent variables to be new variables we reduce the order of a system from $n$ to $n-1$ and increase its dimension form $m$ to $2m$.
Point transformations for the lower order system correspond to generalized contact transformations for the higher order system. We obtained the canonical form of scalar third order ODEs linearizable via generalized contact transformations. This canonical form gave us four equivalence classes for scalar third order ODEs depending on the number of infinitesimal generators. As such, a scalar third order linear ODE has three classes with $4$, $5$ and $10$ contact symmetries. The maximal symmetry class with $15$ generalized contact symmetries corresponds to the maximal symmetry class of contact symmetries. Also of great interest is the study of correspondence between the other classes of generalized contact symmetries and those of contact symmetries. Here we obtained group classification of a scalar third order ODE by reducing it to a system of two second order ODEs. If the reduced system of ODEs is linearizable then it can be solved by using geometric linearization \cite{qadirgeo}.

We can reduce a scalar fourth order ODE to a system of two third order ODEs and following the same procedure of group classification we can find generalized contact symmetries of the scalar ODE. Similarly we can take it to higher order ODEs to get equivalence classes of these equations by simply reducing scalar ODEs to systems of ODEs. One could use the given procedure to find the equivalence classes of systems of higher order ODEs. We can reduce a scalar fourth order ODE to a system of two second order ODEs by double reduction. Similarly any system of ODEs of order $n\geq 3$ can be reduced in steps to a system of second order ODEs to use the benefit of geometric linearization. This procedure can also be employed in steps to relate higher order transformations to point transformations \cite{dutt1}. There is much work that needs to be done in this direction. The equations linearizble via these transformations may form a new class of ODEs that do not fall into IM or IMS or Meleshko classes of linearizable ODEs.

\section*{Acknowledgments}
We are grateful to the Higher Education
Commission (HEC) of Pakistan for support under their project no. 3054.

\end{document}